\documentclass{article}
\usepackage{spconf,amsmath,graphicx}
\graphicspath{{./img/}}

\usepackage{amssymb} 
\usepackage{amsfonts}
\usepackage{amsthm}
\usepackage{booktabs}
\usepackage{url}

\newcommand{\ie}{{\it i.e.}}

\newcommand{\reals}{{\mathbb R}}

\newcommand{\minimize}{\operatornamewithlimits{minimize}}

\title{Robust dissimilarity measure for Network Localization}
\name{Cl\'{a}udia Soares, Jo\~{a}o Gomes \thanks{This research was
    supported by Funda\c{c}\~{a}o para a Ci\^{e}ncia e a Tecnologia
(projects PEst-OE/EEI/LA0009/2013 and PhD grant SFRH/BD/72521/2010) and EU FP7 project MORPH (grant agreement no.\ 288704).}}
\address{Institute for Systems
    and Robotics (ISR), Instituto Superior T\'ecnico, Universidade de
    Lisboa, \\ 1049-001 Lisboa, Portugal\\\texttt{\{csoares,jpg\}@isr.ist.utl.pt}}
\begin{document}
%
\maketitle
\begin{abstract}
  In practice, network applications have to deal with failing nodes,
   malicious attacks, or, somehow, nodes facing highly
  corrupted data --- generally classified as outliers. This calls for
  robust, uncomplicated, and efficient methods. We propose a
  dissimilarity model for network localization which is robust to
  high-power noise, but also discriminative in the presence of regular
  gaussian noise. We capitalize on the known properties of the
  M-estimator Huber penalty function to obtain a robust, but nonconvex,
  problem, and devise a convex underestimator, tight in the function
  terms, that can be minimized in polynomial time. Simulations show
  the performance advantage of using this dissimilarity model in the
  presence of outliers and under
   regular gaussian noise.
\end{abstract}
\section{Introduction}
\label{sec:intro}
Low signal to noise ratio, multipath components, reflections,
interference: these are common hurdles to overcome throughout many
applications of sensor networks --- and one of their major
consequences is the presence of outliers in the collected
data. Outliers are measurements with unexpected and surprising values
given the overall behavior of the sensor network. They can cause large
estimation errors when algorithms are not prepared to account for
them, and if these estimates are an input to some other problem, then
there is a risk that error propagation will invalidate the final
purpose. This is the case with network localization; It might be taken
for granted in most sensor network applications, but it is still a
very active research field. We present a new approach that addresses
the presence of outliers, not by eliminating them of the estimation
process, but by appropriately weighting them, so that they can
contribute to the solution, while mitigating the bias of the
estimator.

\subsection{Related work}
\label{sec:related-work}

Some approaches to robust localization rely on identifying outliers
from regular data. Then, outliers are removed from the estimation of
sensor positions. The work in~\cite{IhlerFisherMosesWillsky2005}
formulates the network localization problem as an inference problem in
a graphical model. To approximate an outlier process the authors add a
high-variance gaussian to the gaussian mixtures and employ
nonparametric belief propagation to approximate the solution. In the
same vein,~\cite{AshMoses2005} employs the EM algorithm to jointly
estimate outliers and sensor positions. Recently, the
work~\cite{YinZoubirFritscheGustafsson2013} tackled robust
localization with estimation of positions, mixture parameters, and
outlier noise model for unknown propagation conditions.

Alternatively, methods may perform a soft rejection of outliers, still
allowing them to contribute to the solution. In the
work~\cite{OguzGomesXavierOliveira2011} a maximum likelihood estimator
for laplacian noise was derived and subsequently relaxed to a convex
program by linearization and dropping a rank constraint, The authors
in~\cite{ForeroGiannakis2012} present a robust Multidimensional
Scaling based on the least-trimmed squares criterion minimizing the
squares of the smallest residuals.  In~\cite{KorkmazVeen2009} the
authors use the Huber loss~\cite{huber1964} composed with a
discrepancy between measurements and estimate distances, in order to
achieve robustness to outliers. The resulting cost is nonconvex, and
optimized by means of the Majorization-Minimization technique.

The cost function we present incorporates outliers into the estimation
process and does not assume any outlier model. We capitalize on the
robust estimation properties of the Huber function but,
unlike~\cite{KorkmazVeen2009}, we do not address the nonconvex cost in
our proposal. Instead, we produce a convex relaxation which
numerically outperforms other natural formulations of the problem.

\subsection{Contributions}
\label{sec:contributions}

We present a tight convex underestimator to each term of the robust
discrepancy measure for sensor network localization. Our approach
assumes no specific outlier model, and all measurements contribute to
the estimate. Numerical simulations illustrate the quality of the
convex underestimator.

\section{The problem}
\label{sec:problem}

The network is represented as an undirected graph~$\mathcal{G} =
(\mathcal{V},\mathcal{E})$. We represent the set of sensors with
unknown positions as~$\mathcal{V} = \{1,2, \dots, n\}$ There is an
edge $i \sim j \in {\mathcal E}$ between nodes $i$ and $j$ if a range
measurement between~$i$ and $j$ is available and~$i$ and $j$ can
communicate with each other.  Anchors have known positions and are
collected in the set ${\mathcal A} = \{ 1, \ldots, m \}$; they are not
nodes on the graph~$\mathcal{G}$. For each sensor $i \in {\mathcal
  V}$, we let ${\mathcal A}_i \subset {\mathcal A}$ be the subset of
anchors with measured range to node~$i$. The set~$N_{i}$ collects the
neighbor sensor nodes of node~$i$.

The element positions belong to~$\reals^p$ with~$p=2$ for planar
networks, and $p=3$ for volumetric ones.  We denote by $x_i \in \reals^p$ the
position of sensor $i$, and by $d_{ij}$ the range measurement between
sensors $i$ and $j$. Anchor positions
are denoted by $a_{k} \in \reals^{p}$. We let $r_{ik}$ denote the
noisy range measurement between sensor $i$ and anchor $k$.

We aim at estimating the sensor positions~$x=\{x_{\mathcal{V}}\}$,
taking into account two types of noise: (1) regular gaussian noise,
and (2) outlier induced noise.

\section{Discrepancy measure}
\label{sec:discrepancy-measure}

The maximum likelihood estimator for the sensor positions with
additive i.i.d.\ gaussian noise contaminating range measurements is the
solution of the optimization problem
\begin{equation*}
  \label{eq:snlOptimizationProblem} 
  \minimize_{x} g_{G}(x),
\end{equation*} 
where
\begin{align*} 
\label{eq:sqdist}
  g_{G}(x) = & \sum _{i \sim j} \frac{1}{2}(\|x_{i} - x_{j}\| -
  d_{ij})^2 + \\
 &\sum_{i} \sum_{k \in \mathcal{A}_{i}} \frac{1}{2}(\|x_{i}-a_{k}\| -
 r_{ik})^2.
\end{align*} 
However, outlier measurements will heavily bias the solutions of the
optimization problem since their magnitude will be amplified by the
squares~$h_{Q}(t) = t^{2}$ at each outlier term.  From robust
estimation, we know some alternatives to perform soft rejection of
outliers, namely, using~$L_{1}$ loss~$h_{|\cdot|}(t) = |t|$ or the
Huber loss
\begin{equation}
\label{eq:huber-loss}
h_{R}(t) =
\begin{cases}
  t^{2} & \text{if } |t| \leq R,\\
  2R|t|-R^{2} & \text{if } |t| \geq R.
\end{cases}
\end{equation}
The Huber loss joins the best of two worlds: it is robust for large values
of the argument --- like the~$L_{1}$ loss --- and for reasonable noise
levels it behaves like~$g_{Q}$, thus leading to the maximum likelihood
estimator adapted to regular noise.
\begin{figure}[tb]
  \centering
  \includegraphics[width=\columnwidth]{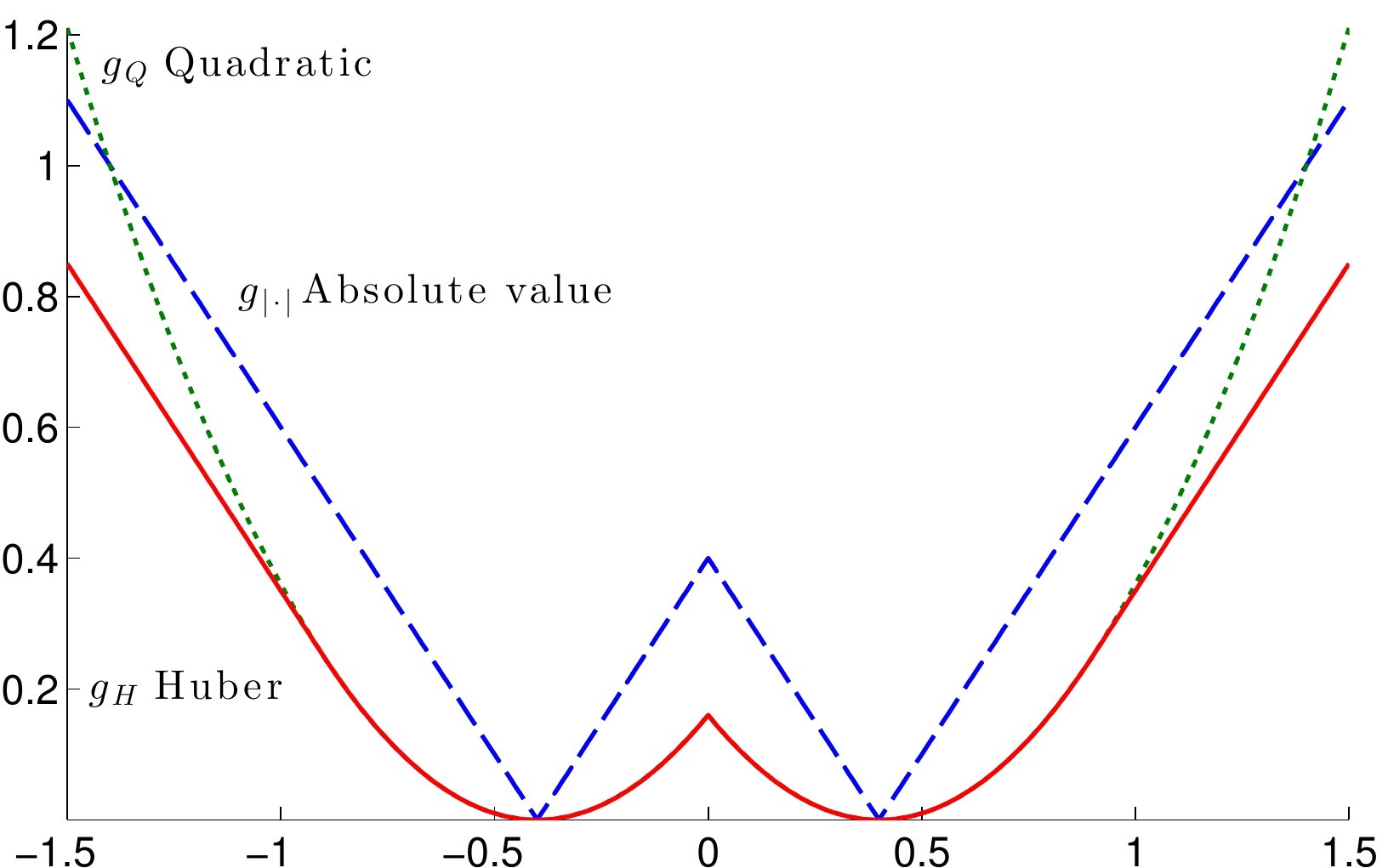}
  \caption{The different cost functions considered in this paper: the
    maximum likelihood independent white gaussian noise~$g_{Q}(x_{i},x_{j}) =
    (\|x_{i}-x_{j}\| - d_{ij})^{2}$ shows the steepest tails,
  which act as outlier amplifiers; the $L_{1}$
  loss~$g_{|\cdot|}(x_{i},x_{j}) = |\|x_{i}-x_{j}\| - d_{ij}|$, associated with
  impulsive noise, which fails to model the gaussianity of regular
  operating noise; and, finally, the Huber loss~$g_{H}(x_{i},x_{j}) = h_{R}(\|x_{i}-x_{j}\| - d_{ij})$, combines robustness
  to high-power outliers and adaptation to medium-power gaussian noise.}
  \label{fig:nonconvex}
\end{figure}
Figure~\ref{fig:nonconvex} depicts a one-dimensional example of these different
costs. We can observe in this simple example the main properties of
the different cost functions, in terms of adaptation to low/medium-power
gaussian noise and high-power outlier spikes.
Using~\eqref{eq:huber-loss} we can write our optimization problem as
\begin{equation}
  \label{eq:snlOptProb}
  \minimize_{x} g_{R}(x)
\end{equation}
where
\begin{align}
  \nonumber
    g_{R}(x) = & \sum _{i \sim j} h_{R_{ij}}(\|x_{i} - x_{j}\| - d_{ij})
    + \\ \label{eq:huberDist}
    &\sum_{i} \sum_{k \in \mathcal{A}_{i}} h_{R_{ik}}(\|x_{i}-a_{k}\| - r_{ik}).
\end{align}
This function is nonconvex and, in general, difficult to minimize.
We shall provide a convex underestimator, that tightly bounds each term
of~\eqref{eq:huberDist}, thus leading to better estimation results
than other relaxations which are not tight~\cite{SimonettoLeus2014}.

\section{Convex underestimator}
\label{sec:conv-under}
To convexify~$g_{R}$ we can replace each term by its
convex hull\footnote{The convex
  hull of a function~$\gamma$, \ie, its best
  possible convex underestimator, is defined as $ \text{conv } \gamma(x) =
  \sup\left \{ \eta(x) \; : \; \eta \leq \gamma, \; \eta \text{ is
      convex} \right \} $. It is hard to determine in general~\cite{UrrutyMarechal1993}.},
\begin{figure}[tb]
  \centering
  \includegraphics[width=\columnwidth]{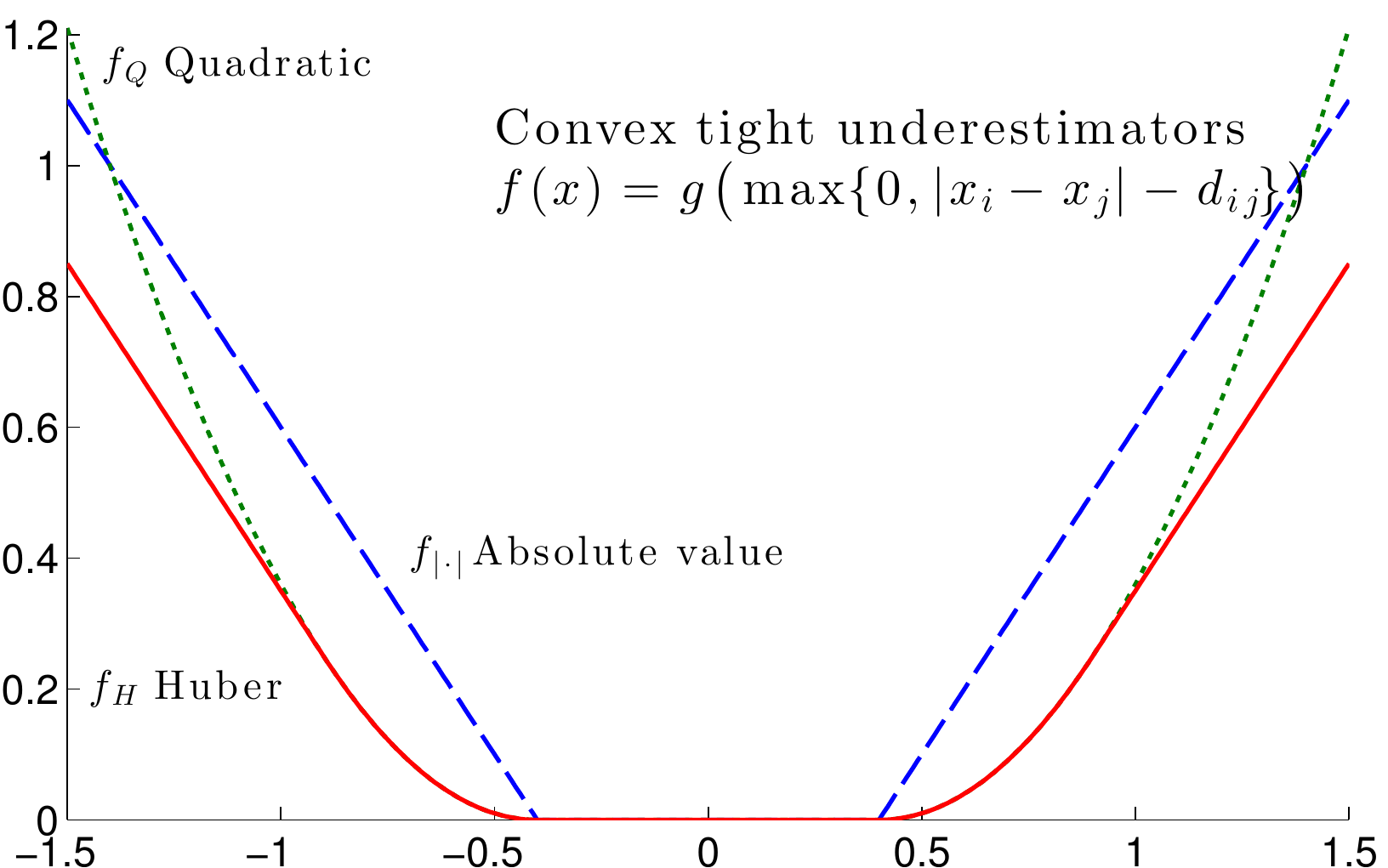}
  \caption{All functions~$f$ are tight underestimators to the
    functions~$g$ in Figure~\ref{fig:nonconvex}. They are the convex
    envelopes and, thus, the best convex approximations to the original
    nonconvex costs. The convexification is performed by restricting
    the arguments of~$g$ to be nonnegative.}
  \label{fig:convexified}
\end{figure}
as depicted in Figure~\ref{fig:convexified}. Here, we observe that the
high-power behavior is maintained, whereas the medium/low-power is
only altered in the convexified area. We define the convex costs by
composition of any of the convex functions~$h$ with a nondecreasing
function~$s$
 \begin{equation*}
   s(t) = \max \{ 0, t\}
 \end{equation*}
which, in turn, transforms the discrepancies
\begin{equation*}
  \delta_{ij}(x) = \|x_{i} - x_{j}\| - d_{ik}, \qquad
  \delta_{ik}(x_{i}) = \|x_{i} - a_{k}\| - r_{ik}.
\end{equation*}
As~$s\left (\delta_{ij}(x)\right )$ and~$s\left (\delta_{ik}(x)\right
)$ are nondecreasing and each one of the functions~$h$ is convex,
then
\begin{align*}
  f(x) = &\sum _{i \sim j} h\left(s\left(\|x_{i} - x_{j}\| -
      d_{ij}\right) \right) + \\
  &\sum_{i} \sum_{k \in \mathcal{A}_{i}} h\left(s\left(\|x_{i}-a_{k}\|
      - r_{ik}\right) \right)
\end{align*}
is also convex. The quality of the convexified quadratic problem was
addressed in~\cite{soares2014simple} and the analysis for the
remaining functions is similar.

\section{Numerical experiments}
\label{sec:numer-exper}

\begin{figure}[tb]
  \centering
  \includegraphics[width=\columnwidth]{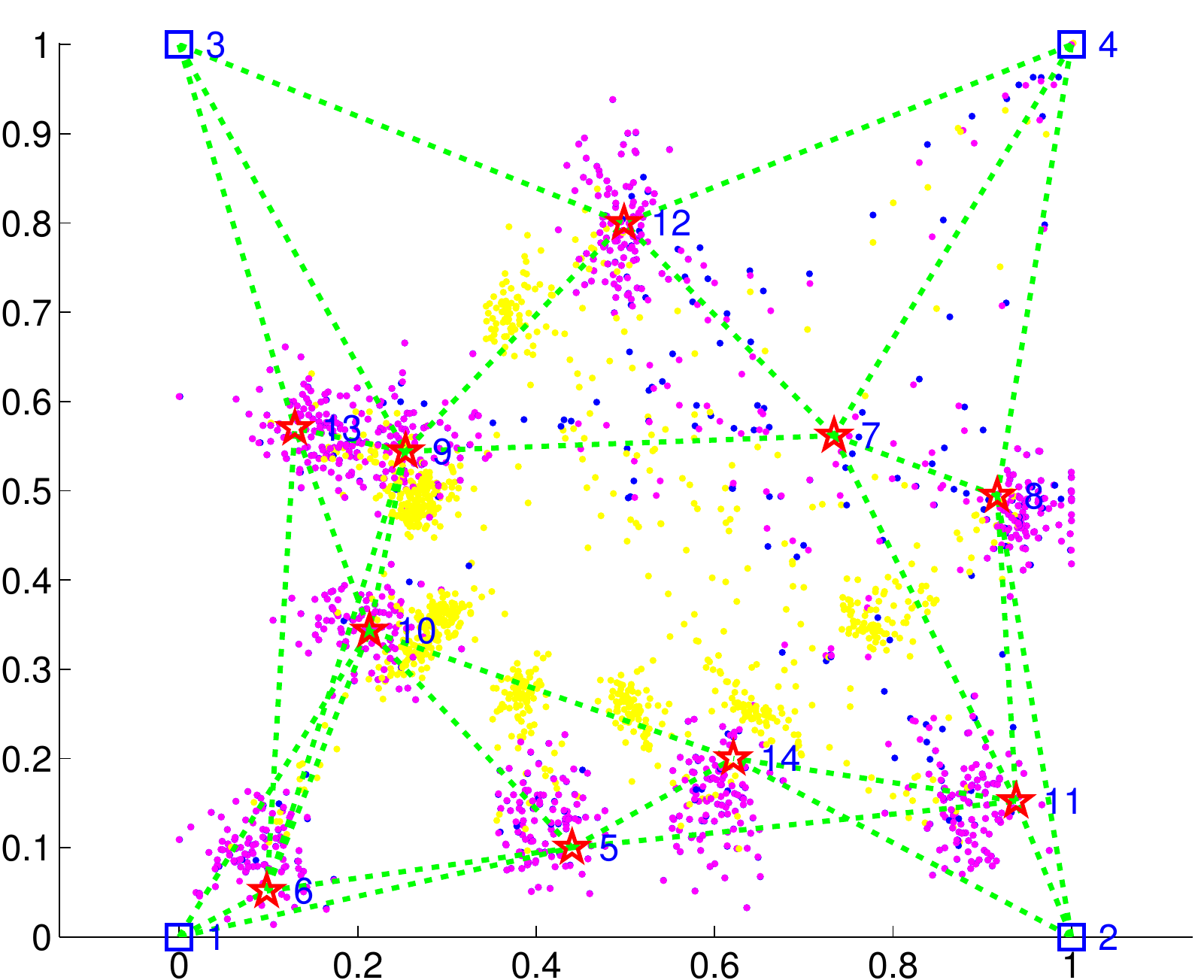}
  \caption{Estimates for sensor positions for the three loss
    functions; We plotted in yellow the monte carlo results of
    minimizing~$g_{|\cdot|}$, the $L_{1}$ loss; in blue we can see the
    estimates resulting from minimizing~$g_{Q}$; in the same way,
    magenta dots represent the output for function~$g_{R}$. It is
    noticeable that the~$L_{1}$ loss is not able to correctly estimate
    positions whose measurements are corrupted with gaussian
    noise. The perturbation in node~$7$ has more impact in the
    augmented dispersion of blue dots than magenta dots around its
    neighbors.}
  \label{fig:estimates}
\end{figure}

We assess the performance of the three considered loss functions
through simulation.  The experimental setup consists in a uniquely
localizable geometric network deployed in a square area with side
of~$1$Km, with four anchors (blue squares in
Figure~\ref{fig:estimates}) located at the corners, and ten sensors,
(red stars). Measurements are also visible as dotted green lines. The
average node degree\footnote{To characterize the network we use the
  concepts of \emph{node degree}~$k_{i}$, which is the number of edges
  connected to node~$i$, and \emph{average node degree}~$<k> = 1/n
  \sum_{i=1}^{n}k_{i}$.} of the network is~$4.3$.  The regular noisy
range measurements are generated according to
\begin{align}
  \label{eq:noise}
  d_{ij} &=& | \|x_{i}^{\star} - x_{j}^{\star}\| + \nu_{ij} |, \; r_{ik} &=& | \|x_{i}^{\star} - a_{k}\| + \nu_{ik} |,
\end{align}
where~$x_{i}^{\star}$ is the true position of node~$i$, and~$\{\nu_{ij} :
i \sim j \in \mathcal{E}\} \cup \{\nu_{ik} : i \in \mathcal{V}, k \in
\mathcal{A}_{i}\}$ are independent gaussian random variables with zero
mean and standard deviation~$0.04$, corresponding to an uncertainty of
about~$40$m.
Node~$7$ is malfunctioning and all measurements related to it are
perturbed with gaussian noise with standard deviation~$4$,
corresponding to an uncertainty of~$4$Km. The convex optimization
problems were solved with \texttt{cvx}~\cite{cvx}.
We ran~$100$ Monte Carlo trials, sampling both regular and outlier
noise.

The performance metric used to assess accuracy is the average
positioning error defined as
\begin{equation}
  \label{eq:avgerror}
  \epsilon = \frac1{MC} \sum_{mc =1}^{MC} \|\hat x(mc) - x^{\star}\|.
\end{equation}
where~$MC$ is the number of Monte Carlo trials performed and~$\hat
x(mc)$ corresponds to the full set of estimates at trial~$mc$.
In Figure~\ref{fig:estimates} we can observe that clouds of estimates
from~$g_{R}$ and~$g_{Q}$ gather around the true positions, except for
the malfunctioning node~$7$. Note the spread of blue dots in the
surroundings of the edges connecting node~$7$, indicating that~$g_{R}$
better preserves the nodes' ability to localize themselves, despite
their confusing neighbor~$7$.
\begin{table}[tb]
  \caption{Average positioning error per sensor~($\epsilon$/sensor), in meters}
  \label{tab:error}
  \centering
  \begin{tabular}[h]{@{}ccc@{}}
    \toprule
    \textbf{$f_{|\cdot|}$} &
    \textbf{$f_{Q}$} &
    \textbf{$f_{R}$}\\\midrule
    59.50&32.16&31.06\\
    \bottomrule
  \end{tabular}
\end{table}
This intuition is confirmed by the analysis of the data in
Table~\ref{tab:error}, which demonstrates that, even with only one
disrupted sensor, our robust cost can reduce the error per sensor
by~$1.1$ meters.
Also, as expected, the malfunctioning node cannot be positioned.
\begin{figure}[tb]
  \centering
  \includegraphics[width=\columnwidth]{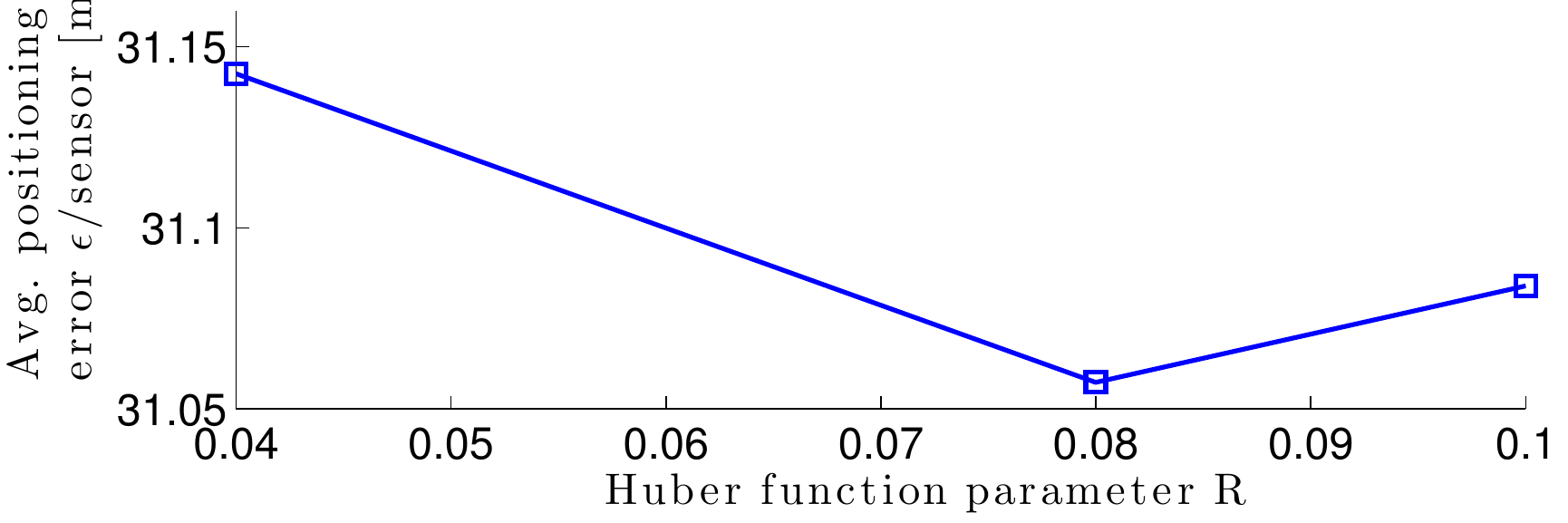}
  \caption{Average positioning error \textit{versus} the value of the Huber
    function parameter~$R$. The accuracy is maintained even for
    different parameter values.}
  \label{fig:varR}
\end{figure}
The sensitivity to the value of the Huber parameter~$R$
in~\eqref{eq:huber-loss} is moderate, as shown in
Figure~\ref{fig:varR}. In fact, the error per sensor of the proposed
estimator is always the smallest for all tested values of the
parameter. We observe that the error increases when~$R$ approaches the
standard deviation of the regular gaussian noise, meaning that the
Huber loss gets closer to the~$L_{1}$ loss and, thus, is no longer
adapted to the regular noise ($R=0$ corresponds exactly to the~$L_{1}$
loss); in the same way, as~$R$ increases, so does the quadratic
section, and the estimator gets less robust to outliers, so, again,
the error increases.

Another interesting experiment is to see what happens when the faulty
sensor produces measurements with consistent errors or bias. So, we
ran~$100$ Monte Carlo trials in the same setting, but node~$7$
measurements are consistently~$10\%$ of the real distance to each
neighbor.
\begin{table}[tb]
  \caption{Average positioning error per sensor~($\epsilon$/sensor),
    in meters, for the biased experiment}
  \label{tab:error-bias}
  \centering
  \begin{tabular}[h]{@{}ccc@{}}
    \toprule
    \textbf{$f_{|\cdot|}$} &
    \textbf{$f_{Q}$} &
    \textbf{$f_{R}$}\\\midrule
    80.98 &  58.31 &  47.08\\ 
    \bottomrule
  \end{tabular}
\end{table}
The average positioning error per sensor is shown in
Table~\ref{tab:error-bias}. Here we observe a significant performance
gap between the alternative costs, and our formulation proves to be, by
far, superior.

\section{Conclusions and future work}
\label{sec:concl-future-work}

We proposed a dissimilarity model easy to motivate and effective,
which accounts for outliers without prescribing a model for outlier
noise. This dissimilarity model was convexified by means of the convex
envelopes of its terms, leading to a problem with a unique minimum
value attainable in polynomial time.

Different types of algorithms can be designed to attack the discrepancy measure
presented in this work, since the function is continuous and
convex. Due to the distributed nature of networks of sensors or, generically,
agents, we aim at investigating a distributed minimization of the
proposed robust loss. There are also several nice properties regarding
distributed operation: the adjustable Huber
parameter which is local to each edge and can be dynamically adjusted
to the local environmental noise conditions, in a distributed manner.

Nevertheless, distributed implementation is a challenge, since the
proposed loss lacks smoothness properties which were valuable in our
previous work~\cite{soares2014simple}.

\bibliographystyle{IEEEbib}
\bibliography{biblos}

\end{document}